\documentclass[twocolum]{autart}

\usepackage{graphics} 
\usepackage{epsfig} 
\usepackage{mathptmx} 
\usepackage{times} 
\usepackage{amsmath} 
\usepackage{amssymb}  
\usepackage{subfigure}
\usepackage{verbatim}
\usepackage{cite}
\usepackage{algorithmic}

\usepackage{savesym}
\savesymbol{AND}

\usepackage{algorithm}
\usepackage{mathtools}

\newcounter{prob1}
\newcounter{prob2}
\newcounter{prob3}
\newcounter{prob4}
\newcounter{prob5}
\newcounter{prob6}
\newcounter{prob7}

\newcommand{\Proof}{{\bf Proof. }}
\newtheorem{ass}{Assumption}

\date{}

\begin{document}
\begin{frontmatter}
\title{Stochastic Sensor Scheduling via Distributed Convex Optimization
\thanksref{footnoteinfo}}
\thanks[footnoteinfo]{This research has been supported partially supported by NSF ECS-0901846 and NSF CCF-1320643. Partial version of this paper has appeared in \cite{Li2011}.}

\author[qual]{Chong Li}\ead{chongli@iastate.edu},
\author[isu]{Nicola Elia}\ead{nelia@iastate.edu}
\address[qual]{Dept. of Electrical and
Computer Engineering, Iowa State University, Ames, IA, 50011. Currently with Qualcomm Research, Bridgewater, NJ, 08807.}
\address[isu]{Dept. of Electrical and Computer Engineering, Iowa State University, Ames, IA, 50011.}


\begin{abstract}
In this paper, we propose a stochastic scheduling strategy for estimating the states of $N$ discrete-time linear time invariant (DTLTI) dynamic systems, where only one system can be observed by the sensor at each time instant due to practical resource constraints.
The idea of our stochastic strategy is that a system is randomly selected for observation at each time instant according to a pre-assigned probability distribution.
We aim to find the optimal pre-assigned probability in order to minimize the maximal estimate error covariance among dynamic systems.
We first show that under mild conditions, the stochastic scheduling problem gives an upper bound on the performance of the optimal sensor selection problem, notoriously difficult to solve.
We next relax the stochastic scheduling problem into a tractable suboptimal quasi-convex form. We then show that
the new problem can be decomposed into coupled small convex optimization problems, and it
can be solved in a distributed fashion.
Finally, for scheduling implementation, we propose centralized and distributed deterministic scheduling strategies
based on the optimal stochastic solution and provide simulation examples.

\end{abstract}
\end{frontmatter}

\setcounter{prob1}{1}
\setcounter{prob2}{2}
\setcounter{prob3}{3}
\setcounter{prob4}{4}
\setcounter{prob5}{5}
\setcounter{prob6}{6}
\setcounter{prob7}{7}

\section{Introduction}
 In this paper, we consider the problem of scheduling the observations of independent targets in order to minimize the
tracking error covariance, but when only one target can be observed at a given time. This problem captures many interesting
tracking/estimation application problems. As a motivational example, consider $N$ independent dynamic targets, spatially distributed in an area,
that need to be tracked (estimated) by a single (mobile) camera sensor. The camera has limited
sensing range and therefore it needs to zoom in on, or be in proximity of, one of the targets for obtaining measurements.
Under the assumption that the switching time among the targets is negligible,  then we need to find a visiting sequence in order to minimize the estimate error.\\
Another case is when a set of $N$ mobile surveillance devices need to track $N$ geographically-separated targets,
where each target is tracked by one assigned surveillance device. However, the sensing/measuring channel can only be used by one
estimator at the time (e.g. sonar range-finding \cite{Murray02}). Then, we need to design a scheduling sequence of surveillance
devices for accurate tracking.

\subsection{Related Work and Contributions of This Paper}
 There has been considerable research effort devoted to the study of sensor selection problems, including sensor scheduling \cite{A.Tiwari04, Alriksson05, Vasanthi06,Liang07,gupta06,Shi07,Birar09,boyd09,Dahleh11, srivastava11randomized, Lin13,zhenliang_Submodularity,zhenliang_Near_optimality,zhenliang_String} and sensor coverage \cite{H.Choset01,E.U.Acar02,J.Cortes04,Cortes10, gupta06, I.I.Hussein07}. This trend has been inspired by the significance and wide applications of sensor networks. As the literature is vast, we list a few results which are relevant to this paper. The sensor scheduling problem mainly arises from minimization of two relevant costs: sensor network energy consumption and estimate error. On the one hand, \cite{Vasanthi06},\cite{Liang07} and \cite{Birar09}, see also reference therein, have proposed various efficient sensor scheduling algorithms to minimize the sensor network energy consumption and consequently maximize the network lifetime. On the other hand, researchers have proposed many tree-search based sensor scheduling algorithms (mostly in conjunction with Kalman filtering) to minimize the estimate error \cite{A.Tiwari04},\cite{Alriksson05},\cite{Lin13}, e.g. sliding-window, thresholding, relaxed dynamic programming, etc. By taking both sensor network lifetime and estimate accuracy into account, several sensor tree-search based scheduling algorithms have been proposed in \cite{Shi07}, \cite{Shi08}. In \cite{boyd09}, the authors have formulated the general sensor selection problem and solved it by relaxing it to a convex optimization problem. The general formulation therein can handle various performance criteria and topology constraints. However, the framework in \cite{boyd09} is only suitable for static systems instead of dynamic systems which are mostly considered in the literature.\\
  In general, deterministic optimal sensor selection problems are notoriously difficult. In this paper, we propose a stochastic scheduling strategy.
At each time instant, a target is randomly chosen to be measured according to a pre-assigned probability distribution.
We find the optimal pre-assigned distribution that  minimize an upper bound on the expected estimate error covariance (in the limit)
in order to keep the actual estimate error covariance small. Compared with algorithms in the literature, this strategy has low computational complexity,
it is simple to implement and provides performance guarantee on the general deterministic scheduling problem. Of course, the reduction of computational complexity comes at the
expenses of degradation of the ideal performance. However, in many situations the extra computational complexity cost may not be justified. Further this strategy can easily incorporate extra constraints on the
scheduling design, which might be difficult to handle in existing algorithms (e.g. tree-search based algorithms).\\
Our work is related to \cite{gupta06}, \cite{Mo_stochastic11} and \cite{Dahleh11}.
\cite{gupta06} introduces stochastic scheduling to deal with sensor selection and coverage problems, and \cite{Mo_stochastic11} extends
the setting and results in \cite{gupta06} to a tree topology.
Although we also adopt the stochastic scheduling approach, the problem formulation and proposed algorithms of this paper
are different from \cite{gupta06,Mo_stochastic11}.
In particular,  we consider different cost functions and design distributed algorithms that provide optimal probability
distributions.
\cite{Dahleh11} has considered a scheduling problem in continuous-time and proposed a tractable relaxation, which provides a lower bound on the achievable performance,
and an open-loop periodic switching strategy to achieve the bound in the limit of arbitrarily fast switching.
However, besides the difference in the formulations, their approach does not appear to be directly extendable to the discrete-time setting. In summary, our main contributions include:
\begin{enumerate}
\item We obtain a stochastic scheduling strategy with performance guarantee on the general deterministic scheduling problem by solving distributed optimization problems.
\item For scheduling implementation, we propose both centralized and distributed deterministic scheduling strategies.
\end{enumerate}

\subsection{Notations and Organization}
 Throughout the paper, $A^{'}$ is the transpose of matrix $A$. $Ones(n,n)$ implies an $n\times n$ matrix with $1$ as all its entries. $Diag(V)$ denotes a diagonal matrix with vector $V$ as its diagonal entries. $M \succeq 0$ (or $M\in S_+$) and $M \succ 0$ (or $M\in S_{++}$) respectively implies matrix $M$ is positive semi-definite and positive definite where $S_{+}$ and $S_{++}$ represent the positive semi-definite and positive definite cones. For a matrix A, if the block entry $A_{ij} = A_{ji}'$, we use $(\cdot)$ in the matrix to present block $A_{ji}$.\\
 The paper is organized as follows. In section 2, we mathematically formulate the stochastic scheduling problem. In section 3, we develop an approach and a distributed computing algorithm to solve the optimization problem. In section 4, we present some further results and the extensions of our model. In section 5, we consider the scheduling implementation problem. At last, we present simulations to support our results.
\section{Sensor Scheduling Problem Setup}
Consider a set of N DTLTI systems (targets) evolving according to the equations
\begin{equation}
x_{i}[k+1]=A_{i}x_{i}[k]+w_{i}[k]   \quad i=1,2,\dotsc N
\label{sys_model01}
\end{equation}
where $x_{i}[k]\in \mathbb{R}^{n_i}$ is the process state vector and $w_{i}[k]\in \mathbb{R}^{n_i}$ is assumed to be an independent Gaussian noise with zero mean and covariance matrix $Q_{i}\succ 0$. The initial state $x_{i}[0]$ is assumed to be an independent Gaussian random variable with zero mean and covariance matrix $\pi_{i}[0]$. In practice, each DTLTI system modeled above may represent the dynamic change of a local environment, the trajectory of a mobile vehicle, the varying states of a manufactory machine, etc. As a result of the sensor's limited range of sensing or the congestion of the sensing channel, at time instant $k$, only one system can be observed as
\begin{equation}
\tilde{y}_i[k]=\xi_{i}[k](C_{i}x_{i}[k]+v_{i}[k])
\label{sys_model02}
\end{equation}
where $\xi_{i}[k]$ is the indicator function indicating whether or not the system $i$ is observed at time instant $k$, and accordingly we have constraint\footnote{If we assume \textit{at most} one target is chosen to be measured at each time instant, then we have $\sum_{i=1}^{N} \xi_i \leq 1$. Without loss of generality, in this paper we consider the case that one out of $N$ targets must be chosen at each time instant.} $\sum_{i=1}^{N}\xi_{i}[k]=1$. $v_i[k]\in \mathbb{R}^{p_i}$ is the measurement noise, which is assumed to be independent Gaussian with zero mean and covariance matrix $R_i\succ 0$.
\begin{ass}
For all $i\in\lbrace 1,2, \cdots, N\rbrace$,  the pair $(A_i, Q_i^{1/2})$ is controllable and the pair $(A_{i},C_{i})$ is detectable.
\label{assumption01}
\end{ass}
Denote $\hat{x}_{i}[k]$ as the estimate at time $k$, obtained by a causal estimator for system $i$, which depends on the past and current observations $\lbrace \tilde{y}_i[j]\rbrace_{j=1}^k$.
We begin by considering problem of minimizing (in the limit) the maximal estimate error. The problem can be formulated mathematically as
\begin{equation}
\begin{split}
\displaystyle\min_{\hat{x}_{i}, \lbrace \xi_{i}[j]\rbrace_{j=1}^{\infty}} & \max_{i}\left(\limsup_{T\rightarrow \infty} \frac{1}{T}\sum_{k=1}^T \mathbb{E}[(x_i[k]-\hat{x}_{i}[k])'(x_i[k]-\hat{x}_{i}[k])]\right)\\
\text{s.t.} \quad& \textit{Equation}:(\ref{sys_model01}),(\ref{sys_model02}), \quad i=1,\ldots, N,\\
\quad & \sum_{i=1}^{N}\xi_{i}[k]=1,\\
\end{split}
\label{Intro_OP}
\end{equation}
\\
As the DTLTI systems are assumed to be evolving independently, then for a fixed $\lbrace\xi_i [j]\rbrace_{j=1}^\infty$ the optimal estimator for minimizing the estimate error covariance of system $i$ ($i =1,2,\cdots, N$) is given by a Kalman filter\footnote{This indicates that $N$ parallel estimators, i.e., Kalman filters, are used for estimating $N$ independent DTLTI systems.} whose process of prediction and update are presented as follows \cite{sinopli}. Firstly we define
\begin{equation*}
\begin{split}
&\hat{x}_{i}[k|k] \triangleq \mathbb{E}[x_i[k]|\lbrace \tilde{y}_i[j]\rbrace_{j=1}^k]\\
&P_{i}[k|k]\triangleq \mathbb{E}[(x_i[k]-\hat{x}_{i}[k|k]) (x_i[k]-\hat{x}_{i}[k|k])'|\lbrace \tilde{y}_i[j]\rbrace_{j=1}^k]\\
&\hat{x}_{i}[k+1|k] \triangleq \mathbb{E}[x_i[k+1]|\lbrace \tilde{y}_i[j]\rbrace_{j=1}^k] \\
&P_{i}[k+1|k]\\
&\triangleq \mathbb{E}[(x_i[k+1]-\hat{x}_{i}[k+1|k])(x_i[k+1]-\hat{x}_{i}[k+1|k])'|\lbrace \tilde{y}_i[j]\rbrace_{j=1}^k]\\
\end{split}
\end{equation*}
Then the Kalman filter evolves as
\begin{equation*}
\begin{split}
&\hat{x}_{i}[k+1|k]=A_{i}\hat{x}_{i}[k|k]\\
&P_{i}[k+1|k]=A_{i}P_{i}[k|k]A'+Q_i\\
&\hat{x}_{i}[k+1|k+1]\\
&=\hat{x}_{i}[k+1|k]+\xi_{i}[k+1]K[k+1](y_{i}[k+1]-C_{i}\hat{x}_{i}[k+1|k])\\
&P_{i}[k+1|k+1]=P_{i}[k+1|k]-\xi_{i}[k+1]K_{i}[k+1]C_{i}P_{i}[k+1|k]\\
\end{split}
\end{equation*}
where $K_{i}[k+1]=P_{i}[k+1|k]C_{i}'(C_{i}P_{i}[k+1|k]C_{i}'+R_{i})^{-1}$ is the Kalman gain matrix. After straightforward derivation, we have the covariance of the estimate error evolving as
\begin{eqnarray}\label{errorCov_evolvingEqu}
A_{i}P_{i}[k]A_{i}^{'}+Q_{i}-\xi_{i}[k]A_{i}P_{i}[k]C_{i}^{'}(C_{i}P_{i}[k]C_{i}^{'}+R_{i})^{-1}C_{i}P_{i}[k]A_{i}^{'}\nonumber\\
=P_{i}[k+1]
\end{eqnarray}
where we use the simplified notation $P_i[k]=P_i[k|k-1]$. Note that the error covariance $P_{i}[k+1]$ is a function of sequence $\lbrace\xi_i [j]\rbrace_{j=1}^k$. Moreover, given $\lbrace\xi_i [j]\rbrace_{j=1}^k$, the evolution of the error covariance $P_{i}$ is independent of the measurement values. Substituting the optimal estimator, the problem (\ref{Intro_OP}) is simplified into the following one.\\

\textit{\bf Deterministic Scheduling Problem:}\\
\begin{equation}
\begin{split}
\mu_d=\displaystyle\min_{\lbrace \xi_{i}[j]\rbrace_{j=1}^{\infty}, i=[1..N]} &\quad \max_{i}\left(\limsup_{T\rightarrow \infty} \frac{1}{T}\sum_{k=1}^T Tr(P_i[k])\right)\\
\text{subject to} &\quad \textit{Equation}:(\ref{errorCov_evolvingEqu}), \quad \sum_{i=1}^{N}\xi_{i}[k]=1,\\
\end{split}
\label{Intro_OP01}
\end{equation}
\\
This problem is notoriously difficult to solve. There are not known optimal alternatives to searching all possible schedules
and then pick the optimal
one by complete comparison. However, the procedure is not computationally tractable in practice.
Motivated by this fact, in what follows, we present a stochastic scheduling strategy with advantages summarized below:
\vspace{-5mm}\begin{enumerate}
\item The stochastic scheduling strategy provides an upper bound on the performance of the deterministic scheduling problem under mild conditions, as proved in Theorem \ref{upperBound} next.
\item The stochastic scheduling problem can be easily relaxed into a convex optimization problem, which can be solved efficiently in a distributed fashion.
\item The relaxed problem provides an open-loop strategy, which  has low computing complexity and is simple to implement.
\item Several practical constraints/considerations can be easily incorporated into the stochastic scheduling formulation, as discussed in Section $4$.
\end{enumerate}
\subsection{Problem Formulation: Stochastic Scheduling Strategy}
First of all, we remove the dependence on time instant $k$ and consider $\xi_{i}$ as an independent and identically distribution (i.i.d) Bernoulli random variable with
\begin{equation}
\xi_{i}[k]=\Big\lbrace\begin{array}{cl}
                 1 & \textrm{with probability}\;\; q_{i} \\
                 0 & \textrm{with probability}\;\; 1-q_{i}
               \end{array} \quad i=1,2,\dotsc N
\label{xi}
\end{equation}
for all $k$, where $q_i$ is the probability that the system $i$ is observed at each time instant. As $\sum_{i=1}^{N} \xi_i = 1$, we have $\sum_{i=1}^{N}q_i=1$. Then the stochastic scheduling strategy is that at each time instant a target (i.e., DTLTI system) is randomly chosen for measurements according to a pre-assigned probability distribution $\lbrace q_i\rbrace_{i=1}^{N}$.\\
\indent Notice that the error covariance $P_i[k+1]$ is random as it depends on the randomly chosen sequence $\lbrace \xi_i[j] \rbrace_{j=1}^k$. Thus we need to evaluate the expected estimate error covariance in order to minimize the actual estimate error. Putting above together, we have the stochastic scheduling problem motivated by (\ref{Intro_OP01}) as follows.\\

\textit{\bf Stochastic Scheduling Problem}
\begin{equation}
\begin{split}
\mu_s=\displaystyle\min_{q_i, i=[1..N]} &\quad \max_{i} Tr\left(\lim_{k\rightarrow \infty}\, \mathbb{E}_{\mathbb{\xi}_{i,k}}[P_i[k]]\right)\\
\text{subject to} &\quad \textit{Equation}: (\ref{errorCov_evolvingEqu}), (\ref{xi}),\\
&\quad \sum_{i=1}^{N}q_{i}=1,\;\; 0 \leq q_i \leq 1,\\
\end{split}
\label{Intro_OP02}
\end{equation}
where the expectation is w.r.t $\{\xi_i[1],\ldots \xi_i[k]\}$, which we denote by  $\mathbb{E}_{\mathbb{\xi}_{i,k}}$.\\

\begin{thm}
If $A_i$ is invertible for $i=1,\ldots,n$, the deterministic scheduling performance $\mu_d$ in (\ref{Intro_OP01}) is \textit{almost surely} upper bounded by the stochastic scheduling performance $\mu_s$ in (\ref{Intro_OP02}).
\label{upperBound}
\end{thm}
\vspace{-.3cm}\Proof
See Appendix.\qed
\begin{rem} While Problem (\ref{Intro_OP01}) provides an important motivation, Problem (\ref{Intro_OP02}) is relevant on its own right and can be applied without the restriction on $A_i$ being invertible.
\end{rem}
\subsection{Relaxation}
Problem (\ref{Intro_OP02}) involves the evolution of
$$
\hspace{-4mm}\begin{array}{l}
\mathbb{E}_{\mathbb{\xi}_{i,k}}[P_{i}[k]]=A_{i}\mathbb{E}_{\mathbb{\xi}_{i,k-1}}[P_{i}[k-1]]A_{i}^{'}+Q_{i}-\\
\, q_{i}\mathbb{E}_{\mathbb{\xi}_{i,k-1}}[A_{i}P_{i}[k-1]C_{i}^{'}(C_{i}P_{i}[k-1]C_{i}^{'}+R_{i})^{-1}C_{i}P_{i}[k-1]A_{i}^{'}]
\end{array}
$$
Unfortunately, the right-hand side of the above expression is not easily computable, as it involves the expectation w.r.t. $\mathbb{\xi}_{i,k-1}$,
of a nonlinear recursive expression of $P_{i}[k-1]$.
However, \cite{sinopli} has nicely shown that $\lim_{k\rightarrow \infty}\mathbb{E}_{\mathbb{\xi}_{i,k}}[P_{i}[k]]$ is upper bounded by the fixed point $X_i \succeq 0$ of the following associated MARE,
\begin{equation}
MARE:\,X_i=A_{i}X_{i}A_{i}^{'}+Q_{i}-q_{i}A_{i}X_{i}C_{i}^{'}(C_{i}X_{i}C_{i}^{'}+R_{i})^{-1}C_{i}X_{i}A_{i}^{'}.\\
\label{eq_riccati}
\end{equation}
This result motivates us to minimize $Tr(X_{i})$ as a means to keep $Tr(\lim_{k\rightarrow \infty}\mathbb{E}_{\mathbb{\xi}_{i,k}}[P_{i}[k]])$ itself small.
Specifically, we consider the following optimization problem, denoted by OP, in the rest of the paper.
\begin{eqnarray}\label{OP1}
\text{OP}:&\displaystyle\min_{q_i,i=[1..N]} \quad \max_{i} Tr(X_i)\\
&\text{subject to:}\;\;  \sum_{i=1}^N q_i=1, \;\; q_i^c < q_i \leq 1, \;\; i=1,2,\dotsc N\nonumber
\end{eqnarray}
where $X_i$ is an implicit function of $q_i$, defined by $X_i=g_{q_i}(X_i)$ and
\begin{equation}
\hspace{-1mm}g_{q_{i}}(X_i)=A_{i}X_{i}A_{i}^{'}+Q_{i}-q_{i}A_{i}X_{i}C_{i}^{'}(C_{i}X_{i}C_{i}^{'}+R_{i})^{-1}C_{i}X_{i}A_{i}^{'}
\end{equation}
\begin{rem} $q_i^c$ is the critical value depending on the unstable eigenvalues of $A_i$, where the fixed point $X_i$ exists if and only if the assigned probability $q_i>q_i^c$. We refer interested readers to \cite{sinopli, Mo} for the details on $q_i^c$. In this paper, we assume $\sum_{i=1}^N q_i^c < 1$, otherwise the above optimization problem has no feasible solution, i.e., the upper bound turns out to be infinity. For stable systems, we always have $q_i^c = 0$.
\label{remark_criticalValue}
\end{rem}
Note that the searching space of the above optimization problem is continuous, and it will be shown that the problem is (quasi)-convex. In addition, based on Theorem \ref{upperBound}, the objective value of the above problem is an upper bound on the performance of the deterministic scheduling problem (\ref{Intro_OP01}). In the rest of the paper, we provide efficient distributed algorithms to obtain the optimal solutions of OP.

\section{Minimization of The Maximal Estimate Error among Targets.}
In this section, we first show that OP can be decoupled into $N$ convex optimization problems, which can be solved separately.
Then by utilizing the classical consensus algorithm, we propose a distributed computing algorithm to obtain the optimal solution of OP. First of all, we recall some results on the MARE in \cite{sinopli}
derived under the assumption $Q\succeq 0$.
\begin{lem}
Fix $q\in \mathbb{R}_{(q^c,1]}$, for any initial condition $X_{0}\succeq0$,
\begin{equation*}
\lim_{k\rightarrow\infty}g_{q}^{(k)}(X_{0})=\lim_{k\rightarrow\infty}\underbrace{g_{q}(g_{q}(\cdots g_{q}(X_{0})))}_{k}=X
\end{equation*}
where $X$ is the unique \textit{positive-semidefinite}  fixed point of the MARE, namely, $X=g_{q}(X)$.\\
\label{lem_convg_MARE}
\end{lem}
\begin{lem}
For a given scalar $q$ and a DTLTI system $(A,C,Q,R)$ as described in (\ref{sys_model01}) and (\ref{sys_model02}), the fixed point $X$ of the MARE presented in the form of (\ref{eq_riccati}) can be obtained by solving the following LMI problem.
\begin{equation}
\begin{split}
\quad \underset{X}{\text{argmax}} &\quad Tr(X) \\
\text{subject to}
&\quad \begin{bmatrix}AXA'-X+Q & \sqrt{q}AXC'\\\sqrt{q}CXA'& CXC'+R \end{bmatrix}\succeq 0\\
&\quad X\succeq0.\\
\end{split}
\end{equation}
\label{lem_MARE_LMI}
\end{lem}
\begin{lem}
\label{lem_MARE}
Assume $X, Q\in S_+$,  $R\in S_{++}$ and $(A,Q^{\frac{1}{2}})$ is controllable. Then the following facts are true.
\begin{enumerate}
\item  $g_{q}(X)\succeq g_{q}(Y)$ if $X\succeq Y$.
\item  $g_{q_{1}}(X)\succeq g_{q_{2}}(X)$ if $q_{1}\leq q_{2}$.
\item  $g_{q}(X)=\Phi_{q}(K,X)$ if $K=-AXC'(CXC'+R)^{-1}$.
\item  $g_{q}(X)=\min_{K}\Phi_{q}(K,X)\preceq\Phi_{q}(K,X),\forall K$.
\item  $g_{q}(\alpha X +(1-\alpha)Y)\succeq \alpha g_q (X)+(1-\alpha) g_q(Y), \forall \alpha\in [0,1]$.
\item  Define the linear operator
\begin{equation*}
\quad L_{q}(Y)=(1-q)(AYA')+qFYF'.
\end{equation*}
Suppose there exists $\bar{Y}\succ 0$ such that $\bar{Y} \succ L_{q}(\bar{Y})$. Then for all $W \succeq 0$
\begin{equation*}
\quad \lim_{k\rightarrow \infty} L_{q}^{(k)}(W)=0
\end{equation*}
\end{enumerate}
where
\begin{equation}
\begin{split}
&\Phi_q(K,X)\\
&=(1-q)(AXA'+Q)+q(A+KC)X(A+KC)'+qQ+qKRK'.
\end{split}
\label{equ_Phi}
\end{equation}
\end{lem}
 Now we prove the monotonicity of the fixed point of the MARE, which will facilitate us to analyze OP.
\begin{defn}(\textit{Matrix-monotonicity})
A function $f$: $\mathbb{R} \rightarrow \mathbb{S}_+$ is matrix-monotonic if for all $x, y\in \text{\bf dom} f$ with $x \leq y$, we have either $f(x)\preceq f(y)$ or $ f(x)\succeq f(y)$ in the positive semidefinite cone $\mathbb{S}_+$.
\end{defn}
\begin{thm}
(\textit{Matrix-monotonicity of the MARE}) For $\forall q\in\mathbb{R}_{(q^c,1]}$, the fixed point of the MARE is matrix-monotonically decreasing w.r.t the scalar $q$.
\label{thm_monot}
\end{thm}
\Proof
 Assume $q^c<q_{1}\leq q_{2}\leq 1$ and $X_1$, $X_2$ satisfying $X_{1}=g_{q_{1}}(X_{1})$ and $X_{2}=g_{q_{2}}(X_{2})$. The existence of the fixed points $X_1$ and $X_2$ is guaranteed according to Lemma \ref{lem_convg_MARE}. We need to show $X_1\succeq X_2$. Since $q_{1}\leq q_{2}$, by using Lemma \ref{lem_MARE}(2) we have
\begin{equation*}
X_{1}=g_{q_{1}}(X_{1})\succeq g_{q_{2}}(X_{1}).
\end{equation*}
By Lemma \ref{lem_MARE}(1), we have
\begin{equation*}
\begin{split}
X_{1}&\succeq g_{q_{2}}(g_{q_{2}}(X_{1}))\\
&\succeq g_{q_{2}}(g_{q_{2}}(g_{q_{2}}(X_{1})))\\
&\succeq \dotsc \\
&\succeq g_{q_{2}}^{(k)}(X_{1})\\
\end{split}
\end{equation*}
By the convergence property of the MARE (i.e., Lemma \ref{lem_convg_MARE}), we have $X_{1}\succeq X_{2}$ by taking $k\rightarrow \infty$.\qed
\begin{rem}
This theorem reveals an important message that, for any two different scalar $q_1$ and $q_2$, the corresponding fixed points can be ordered in the positive semidefinite cone. In other words, for a given system model $(A,C,Q,R)$, the fixed points of the MARE w.r.t variable $q$ are comparable. As we will see, this property is the foundation for deriving algorithms to solve OP.
\label{remark_monotonicity}
\end{rem}
Now, we are ready to analyze and solve OP. For the ease of reading, we occasionally use notation $X_i(q_i)$ to stress that $X_i$ is a function of $q_i$, i.e., $X_i$ is the fixed point of the MARE w.r.t the scalar $q_i$.
\begin{cor}
Problem OP is a quasi-convex optimization problem.
\end{cor}
\Proof
Consider the cost function of OP.
From Theorem \ref{thm_monot},  $Tr(X_i(q_i))$ is monotonically non-increasing in $q_i$, as  the trace function is linear. Therefore, $Tr(X_i(q_i))$ is a quasi-convex function, for each $i=1,\ldots,N$, due to the fact that any monotonic function is quasi-convex. Next, based on the fact that non-negative weighted maximum of quasi-convex functions preserves quasi-convexity, the result follows.\qed\\
\\
Next, we can rewrite problem OP in the following equivalent form
\begin{equation}\label{quasiconvex.eq}
\begin{array}{l}
\displaystyle\min_{q_i,i=[1..N],\,\gamma> 0} \quad \gamma\\
subject\; to:\;Tr(X_i(q_i))\leq\gamma,\\
\;\sum_{i=1}^N q_i=1,\;\; q_i^c < q_i \leq 1,\;\; i=1,2,\dotsc N,
\end{array}
\end{equation}
The problem is in principle solved by bisecting $\gamma$ and checking the feasible set is not empty.
However the constraint set is not in a useful form yet. For any fixed $\gamma$, the feasible set is convex but not easy to work with, given the implicit functions $Tr(X_i(q_i))$.

It is then convenient to consider the following related problem for a given $\gamma > 0$.
\begin{equation}\label{lp.eq}
\begin{split}
\mu(\gamma)=\displaystyle\min_{q_i,i=[1..N]}&\quad \sum_{i=1}^N q_i\\
\text{subject to:} & \quad Tr(X_i(q_i))\leq\gamma,\\
&\quad  q_i^c < q_i \leq 1,\quad  \quad i=1,2,\dotsc N,\\
\end{split}
\end{equation}
\begin{lem}\label{feas.lem} $\gamma>0$ is feasible for Problem (\ref{quasiconvex.eq})
if and only if  Problem (\ref{lp.eq}) is feasible and
$\mu(\gamma)\leq 1$.
\end{lem}
\Proof
Let ${\mathcal S}_\gamma$ denote the set  of $q_i$'s feasible for Problem  (\ref{quasiconvex.eq}).

Assume ${\mathcal S}_\gamma$ is not empty. Then feasible set of Problem (\ref{lp.eq}) is not empty, and
since we know that there are $q_i$'s such that  $\sum_{i=1}^N q_i=1$, then $\mu(\gamma)\leq 1$.
For the other direction, assume that $\mu(\gamma)=1$.  Then ${\mathcal S}_\gamma$ is not empty, and therefore $\gamma$ is feasible for Problem (\ref{quasiconvex.eq}).
If $\mu(\gamma)<1$, then let $\alpha>1$ such that $\alpha\mu(\gamma)=1$, and consider $\tilde{q}_i=\alpha q_i$, for $i=1,\ldots,N$.
Then, $\sum_{i=1}^N \tilde{q}_i=1$, $\tilde{q}_i > q_i^c$, and $Tr(X_i(\tilde{q}_i))\leq Tr(X_i(q_i))\leq \gamma$ because of Theorem \ref{thm_monot}.
Thus, ${\mathcal S}_\gamma$ is not empty. Hence $\gamma$ is feasible for Problem (\ref{quasiconvex.eq}).  \qed\\
\indent The cost of (\ref{lp.eq}) is separable and the constraints are independent for each $i$. Thus,
for any $\gamma>0$, (\ref{lp.eq}) can be solved by minimizing $N$ independent problems, namely:
\begin{equation}
\begin{array}{lcccl}
\mu(\gamma)&=&\displaystyle\sum_{i=1}^N &\displaystyle\min_{q_i}& q_i\\
&&&\text{subject to:} &Tr(X_i(q_i))\leq\gamma, \quad q_i^c < q_i \leq 1
\end{array}
\label{op0_equ}
\end{equation}
It is easy to infer that $\mu(\gamma)$ in (\ref{op0_equ}) is decreasing w.r.t. performance $\gamma$.
Thus, OP can be solved by (\ref{op0_equ}) using bisection on $\gamma$ until  $\mu(\gamma)=1$.

We next concentrate on the  subproblems:
\begin{equation}
\begin{split}
q_i^{opt}(\gamma)=\displaystyle\min_{q_i}&\quad q_i\\
\qquad\text{subject to:} &\quad Tr(X_i(q_i))\leq \gamma, \quad q_i^c < q_i \leq 1 \\
\end{split}
\label{op_equ}
\end{equation}
Based on Theorem \ref{thm_monot}, we see that the optimal solution $q_i^{opt}(\gamma)$ of the problem (\ref{op_equ}) implies the smallest probability required for measuring system $i$ for achieving the pre-assigned estimate performance $\gamma$. If the problem (\ref{op_equ}) is not feasible (e.g. $\gamma$ is too small), we set $q_i^{opt}(\gamma) = 1$.
Next, we show that the optimization problem (\ref{op_equ}) can be reformulated as the iteration of a Linear Matrix Inequality (LMI) feasibility problem. Without abuse of notation, we remove the subscript $i$ since the following results apply to all DTLTI dynamic systems.

\begin{lem}
Assume that $(A,Q^{1/2})$ is controllable and $(A,C)$ is detectable. For any given $q\in (q^c,1]$ and invertible matrices $Q$ and $R$, the following statements are equivalent:
\begin{enumerate}
\item $\exists \bar{X} \in S_{++}$ such that $\bar{X}=g_{q}(\bar{X})$.
\item $\exists K$ and $X \in S_{++}$ such that $X \succ \Phi_q(K,X)$ (defined in (\ref{equ_Phi})).
\item $\exists H$ and $G \in S_{++}$ such that $\Gamma_q(H,G)\succ 0$.
\end{enumerate}
where
\begin{equation}\label{Gamma.eq}
\Gamma_q(G,H)=\begin{bmatrix} G &  \sqrt{1-q}GA &  \sqrt{q}(GA+HC)& \sqrt{q}H &G \\ (\cdot)'& G&0&0&0 \\(\cdot)'&(\cdot)'&G&0&0 \\(\cdot)'&(\cdot)'&(\cdot)'&R^{-1}&0\\(\cdot)'&(\cdot)'&(\cdot)'&(\cdot)'&Q^{-1} \end{bmatrix}
\end{equation}
\label{lemma01}
\end{lem}
\vspace{-.8cm}\Proof
$1)\Rightarrow 2)$. According to Lemma \ref{lem_MARE} (3), we have $\bar{X}=g_{q}(\bar{X}) = \Phi_{q}(K_{\bar{X}},\bar{X})$ or $\bar{X}= (1-q)A\bar{X}A'+q(A+KC)\bar{X}(A+KC)'+Q+qKRK'$
with $K_{\bar{X}}=-A\bar{X}C'(C\bar{X}C'+R)^{-1}$.
Then $\bar{X}\succ 0$ since $Q\succ 0$ and the other terms are $\succeq 0$. Moreover,
\begin{equation}
\begin{split}
2\bar{X} =& 2\Phi_{q}(K_{\bar{X}},\bar{X})\\
=& (1-q)A(2\bar{X})A'+q(A+K_{\bar{X}}C)(2\bar{X})(A+K_{\bar{X}}C)'\\
&+2Q+2qK_{\bar{X}}RK_{\bar{X}}'\\
\succ & (1-q)A(2\bar{X})A'+q(A+K_{\bar{X}}C)(2\bar{X})(A+K_{\bar{X}}C)'\\
&+Q+qK_{\bar{X}}RK_{\bar{X}}'\\
=& \Phi_{q}(K_{\bar{X}},2\bar{X})\\
\end{split}
\end{equation}
The inequality follows from the fact that $Q\succ 0$ and $K_{\bar{X}}RK_{\bar{X}}'\succeq 0$. 
The proof is complete.\\
\indent $2)\Rightarrow 1)$. If $X\succ \Phi_q(K,X)$, the proof follows from Theorem $1$ in \cite{sinopli} with $Q\succ 0$. \\
\indent $2)\Leftrightarrow 3)$.
\begin{equation*}
\begin{array}{ll}&X\succ \Phi_q(K,X)\\
\Leftrightarrow &X\succ (1-q)AXA'+q(A+KC)X(A+KC)'+qKRK'+Q
\end{array}
\end{equation*}
Let $G=X^{-1}\succ 0$ and $H=X^{-1}K$. Left and right multiply the above inequality by $G$ we have
\begin{equation*}
\begin{array}{ll}
\Leftrightarrow & GXG\succ (1-q)GAXA'G+q(GA+HC)X(GA+HC)'\\
&+qHRH'+GQG
\end{array}
\end{equation*}
By using Schur complement  this is equivalent to (\ref{Gamma.eq}).\qed
\begin{thm}
If $(A,Q^{1/2})$ is controllable and $(A,C)$ is detectable, the solution of the optimization problem (\ref{op_equ})
can be obtained by solving the following quasi-convex optimization problem in the variables $(q,G,H,Y)$.
\begin{equation}
\begin{split}
\displaystyle\min_{q,H,G\succ 0,Y\succ 0} &\quad q \\
\text{subject to} &\quad Tr(Y)\leq \gamma, \quad i=1,2,\dotsc N \\
&\quad \begin{bmatrix}Y&I\\I&G\end{bmatrix}\succeq 0\\
&\quad \Gamma_{q}(G,H)\succ 0, \quad q^c < q \leq 1,
\end{split}
\label{op_eq01}
\end{equation}
where  $\Gamma_{q}(G,H)$ is given by (\ref{Gamma.eq}).
\end{thm}
\Proof
From the substitution of $G=X^{-1}$ in Lemma \ref{lemma01}, it is straightforward to obtain the following equivalent statements by the Schur complement.\\
$a)$ $\exists X \in S_{++}$ such that $Tr(X)\leq \gamma$.\\
$b)$ $\exists Y \in S_{++}$ such that $Y-X\succeq 0$ and $Tr(Y)\leq \gamma$\\
$c)$ $\exists Y,G \in S_{++}$ such that $Tr(Y)\leq \gamma$ and $
\begin{bmatrix}Y&I\\I&G\end{bmatrix}\succeq 0.$\\
\indent From Lemma \ref{lemma01}, we have the equivalence between $X=g_{q}(X)$ and $\Gamma_{q}(G,H)\succ 0$ in terms of feasibility. For fixed $q$, $\Gamma_{q}(G,H)\succ 0$ is a LMI in variables $(G,H)$. Therefore, the problem (\ref{op_eq01}) can be solved as a quasi-convex optimization problem by using bisection for variable $q$.\qed

In summary, we have shown  the following
\begin{thm}
Problem OP is equivalent to a quasi-convex optimization problem that can be solved by solving (\ref{op0_equ}) using bisection on $\gamma$ until  $\mu(\gamma)=1$ within the desired accuracy.
For  each level  $\gamma$, each of the $N$ independent subproblems  (\ref{op_equ}) can be solved by solving (\ref{op_eq01}) also using bisection for variable $q_i$ and iterating LMI feasibility problems.
\end{thm}

\subsection{Distributed solutions}
Note that the steps of the  outer bisection iteration  are straightforward and can be done either by a centralized scheduler or in a distributed fashion.
If a centralized scheduler/computing-unit is available, it can collect the $q_i^{opt}(\gamma)'s$ from the estimators, check that their sums is less than or equal to one, and send back to the estimator an updated value of $\gamma$ based on a bisection algorithm.

Alternatively, the estimators need to cooperate and agree on an optimal feasible $\gamma$. This can be done assuming the estimators
are strongly connected via a network where the communications between
any two estimators are error-free\footnote{Note that the decentralized computing units are allowed to be allocated in a single fusion center or to be physically distributed
in an area.}. In this case, each estimator needs to obtain the value of $\mu(\gamma)=\sum_{i=1}^N q_i^{opt}(\gamma)$ by communicating with its neighbors.
Under the assumption that $N$ is known to the estimators, $\frac{1}{N}\sum_{i=1}^N q_i^{opt}(\gamma)$ can be obtained by a
distributed averaging process in finite steps as shown in \cite{Sundaram}.
Then, by increasing or decreasing $\gamma$ under a common bisection rule among estimators, the value of $\sum_{i=1}^N q_{i}^{opt}(\gamma)$ can be driven to $1$ and consequently OP is solved.

The above argument, leads to the following distributed computing algorithm, Algorithm \ref{alg_bisection_distributed}, to solve OP.
The inputs of the algorithm are global information assumed to be known by each estimator in prior. Denote $\gamma^{opt}$ as the optimal objective value of OP.
To avoid cumbersome details and to save space, we present the algorithm under the assumption that the interval $[l, u]$ is selected to contain $\gamma^{opt}$.
 i.e., we have $l \leq \gamma^{opt} \leq u$ at each step.
Then the algorithm is guaranteed to converge to the optimal objective value $\gamma^{opt}$ within the desired tolerance.
\begin{algorithm}
\caption{Distributed algorithm for solving OP}
\begin{algorithmic}[1]
\REQUIRE{$N, l\leq \gamma^{opt},  u\geq \gamma^{opt},  \text{tolerance}\quad \epsilon\geq 0$.}
\ENSURE{$\lbrace q_i^{opt}\rbrace_{i=1}^N$, $J^*$}
\FORALL{$i \in\lbrace 1,2,\cdots, N\rbrace$}
\STATE{$u_i \leftarrow u$ and $l_i \leftarrow l$.}
\ENDFOR
\FORALL{$i \in\lbrace 1,2,\cdots, N\rbrace$}
\WHILE{$u_i-l_i > \epsilon$ \COMMENT{Operations in this loop are synchronized among estimators.}}
\STATE{$\gamma \leftarrow \frac{l_i+u_i}{2}$.}
\STATE{Obtain $q_{i}^{opt}(\gamma)$ by solving problem (\ref{op_eq01}).}
\STATE{Obtain $\mu(\gamma) = \sum_{i}q_{i}^{opt}$ via the distributed averaging algorithm.}
\IF{$\mu(\gamma)\leq 1$}
\STATE{$u_i \leftarrow \gamma $}
\ELSE
\STATE{$l_i \leftarrow \gamma $}
\ENDIF
\ENDWHILE
\ENDFOR
\STATE{Objective value of OP $J^* \leftarrow \gamma$}
\end{algorithmic}
\label{alg_bisection_distributed}
\end{algorithm}

\section{Extensions and Special Cases.}\label{extension.sec}
In practical scenarios, some conditions/constraints might be of interest in the sensor scheduling problem.
However, adding extra constraints on scheduling design is problematic in many existing scheduling strategies.
Our stochastic scheduling approach can easily incorporate extra conditions/constaints, as shown through two specific examples next.      \subsection{Prioritization of Certain Targets}
\indent For some reason, specific targets may require extra attention, (i.e., more precise estimation) from the sensor.
Our model can incorporate this requirement by adding constraint $q_j\geq \alpha_j$ where $\alpha_j$ represents the assigned attention weight for target $j$. Then we need to solve the following problem,
\begin{equation}
\begin{split}
\quad \displaystyle\min_{q_i, i=[1..N]} &\quad \max_{i} Tr(X_i) \\
\text{subject to} &\quad X_i=g_{q_i}(X_i),\quad \sum_{i=1}^N q_i=1, \quad q_j\geq \alpha_j\\
&\quad  q_i^c < q_i \leq 1 \quad i=1,2,\dotsc N.\\
\end{split}
\label{ext_op1}
\end{equation}
\subsection{Measurement Loss in Sensing}
\indent In practice, measurement loss is a common phenomena due to various sources, e.g., shadowing, weather condition, large delay, etc. If the measurement loss probability $\tau_i$ for sensing $i$-th target is known in prior, this extra condition can be easily incorporated in our model. Assume that $q_i$'s are pre-assigned to each target. Then the actual probability of reliably receiving measurements from $i$-th target is $q_i(1-\tau_i)$ because of measurement loss. Therefore, OP can be modified as follows,
\begin{equation}
\begin{split}
\quad \displaystyle\min_{q_i, i=[1..N]} &\quad \max_{i} Tr(X_i)\\
\text{subject to} &\quad X_i=g_{q_i(1-\tau_i)}(X_i) \\
&\quad \sum_{i=1}^N q_i=1\\
&\quad q_i(1-\tau_i) >q_i^c, \quad q_i \leq 1, \quad i=1,2,\dotsc N\\
\end{split}
\label{ex_op2}
\end{equation}
With simple modifications, these two extended problems can be solved by the proposed distributed algorithms as well.

\subsection{Closed-form Solutions to MARE for Special Cases}
\indent For the following special class of systems, the underlying MARE has a closed-form solution. Although
this cannot be expected in general, it increases the computational efficiency of the proposed algorithm.  Incidentally,
this appears to the first non-trivial closed-form solution of the MARE in the literature.

\indent Consider a set of $N$ DTLTI single-state systems to be measured evolving according to the equation
\begin{equation}
\begin{split}
x_{i}[k+1]=&a_{i}x_{i}[k]+w_{i}[k]\\
\end{split}
\label{Model_sec4_01}
\end{equation}
where $x_{i}[k],v_i[k],w_{i}[k]\in \mathbb{R}$ and the covariance of $w$ and $v$ are $Q_{i}\in \mathbb{R}^+$ and $R_i\in \mathbb{R}^+$, respectively. The measurement taken by the sensor at each time instant is formulated as follows,
\begin{equation*}
\tilde{y}_{i}[k]=\xi_{i}[k](x_{i}[k-d_i]+v_{i}[k])
\end{equation*}
where $d_{i}$ represents the delay in measurement, which we assume to be fixed and known in this paper. By using augmented states to deal with delays, it is straightforward to have the following compact form for system $i$ with measurement delays,
\begin{equation}
\begin{split}
X_{i}[k+1]&=A_{i}X_{i}[k]+B_{i}w_{i}[k] \\
\tilde{y}_{i}[k]&=\xi_{i}[k](C_{i}X_{i}[k]+v_{i}[k]),\\
\end{split}
\label{Model_sec4_02}
\end{equation}
where $X,A,B,C$ has the following structure\\
\begin{equation}
X_{i}[k]=\begin{bmatrix} x_{i}^{1}[k] \\  x_{i}^{2}[k] \\ \vdots \\  x_{i}^{d_{i}}[k]\\x_{i}[k] \end{bmatrix}, \quad A_{i} = \begin{bmatrix} 0 & 1 & 0 & \cdots & 0\\ 0 & 0 & 1 & \cdots & 0 \\ & \vdots & & \ddots & \vdots \\ 0 & 0 & 0 & \cdots & 1\\ 0 & 0 & 0 & \cdots & a_{i} \end{bmatrix},\quad B_{i} = \begin{bmatrix}0\\0\\ \vdots\\0\\1 \end{bmatrix}
\label{ABC_struc}
\end{equation}
\begin{equation*}
C_{i} = \begin{bmatrix}1&0&\cdots&0\end{bmatrix}.
\end{equation*}
\indent Note that only $x_{i}[k]$ is the true state of system $i$ at time instant $k$ while other states included in vector $X_i[k]$ are dummy variables for handling delays. By exploiting the special structure of above model, we are able to obtain the closed-form fixed point of the MARE. We present this result in the following theorem.
\begin{thm}
For a given $ 0 \leq q \leq 1$, consider the MARE described as $(A,C,R,\tilde{Q})$, where $(A,C)$ have the structure presented in (\ref{ABC_struc}), $R\in \mathbb{R}^+$ and $\tilde{Q} = BQB'$ with $Q\in \mathbb{R}^+$. Then the MARE has a unique positive-semidefinite fixed point $X$ as follows,
if $a=1$ and $ q \neq 0$,
\begin{equation}
X=\begin{bmatrix}x_{1}&x_{1}& \cdots &x_{1} \\ x_{1} & x_{2}& \cdots & x_{2}\\ \vdots & \vdots & \ddots &\vdots& \\
x_{1}& x_{2}& \cdots & x_{n}\end{bmatrix}
\label{sol_MARE01}
\end{equation}
where
\begin{equation*}
x_{j}=\frac{Q+\sqrt{Q^{2}+4qQR}}{2q}+(j-1)Q \quad    j=1,2,\cdots,n ;
\end{equation*}
if  $0< a < \sqrt{\frac{1}{1-q}}$ and $a \neq 1$,
\begin{equation}
X=\begin{bmatrix}x_{1}&ax_{1}& \cdots &a^{n-1}x_{1} \\ax_{1}& x_{2}& \cdots & a^{n-2}x_{2}\\\vdots&\vdots &\ddots &\vdots& \\
a^{n-1}x_{1} &a^{n-2}x_{2}& \cdots & x_{n}\end{bmatrix}
\label{sol_MARE02}
\end{equation}
where
\begin{equation*}
\begin{split}
&x_{1}=\frac{Ra^{2}-R+Q+\sqrt{(Ra^{2}-R+Q)^{2}-4(a^{2}-1-a^{2}q)QR}}{2(1+a^{2}q-a^{2})}\\
&x_{j}=a^{2(j-1)}x_{1}+\frac{1-a^{2(j-1)}}{1-a^2}Q  \quad  j=1,2,\dotsc,n; \\
\end{split}
\end{equation*}
if $a \geq \sqrt{\frac{1}{1-q}}$, MARE fails to converge to a steady state value.\\
\label{theorem02}
\end{thm}
The proof is tedious but straightforward by plugging in the above closed-form solution into the MARE.
If $a \geq \sqrt{\frac{1}{1-q}}$, the results directly follows from \cite{Mo}.

Note that the stochastic upper bound in Theorem 1 may not hold since $A$ is
not invertible.
However, a deterministic scheduling sequence can be  constructed based on the
optimal stochastic solution in the next section, providing an upper bound on the original deterministic
scheduling problem with good performance.

\section{Scheduling Implementation}
For completeness, in this section we consider the problem of scheduling implementation and present a simple approach to implement the scheduling sequence.
For stochastic scheduling implementation, a central scheduler is required to construct a scheduling sequence by randomly selecting targets (via a random seed) according to the optimal probability distribution.
Note that this construction process can be performed efficiently either off-line or on-line.

We next turn our attention back to a deterministic scheduling and look for one  consistent with the the optimal stochastic solution.
 Note that the optimization problems are minimizing the average costs of all possible stochastic sequences.
With the optimal stochastic solutions, we are able to randomly construct a sequence compatible with the distribution.
However, in practice, such random-constructed scheduling sequence may result in undesirable performance. For example,
one target may not be measured for a long consecutive time instants and its error covariance is temporarily built up.
Thus, we would like to identify and use, among all possible stochastic sequences, those that have low costs.
Motivated by the sensor scheduling literature, which suggests periodic solutions \cite{gupta_smartSensor07, Zhang_OptSolSensor10}, we define and look for
deterministic sequences of minimal consecutiveness, defined next. These sequences are periodic and switch among targets most
often compatibly with the optimal scheduling distribution.
We remark that the approach we propose in this section is heuristic but it can be implemented in a distributed fashion and leads to good performance in simulations. We leave the analysis of this and other approaches to future research.

\begin{defn}
Let  $\lbrace s[k]\rbrace_{k=1}^{L}$ be a set of sequences with length $L$, where each element $s[k]$ in the sequence takes value from an element set $\mathbb{K} = \lbrace a_1, a_2, \cdots, a_N \rbrace$ and
the number of occurrences of each value $a_i$ in the sequence is $n_i$. Under these assumptions, then the sequence of minimal consecutiveness is the solution of the following optimization problem
\begin{equation*}
\min_{\lbrace s[k]\rbrace_{k=1}^{L}}\max_{i,j\in\lbrace 1,2,\cdots,N\rbrace} \lbrace j-i|j\geq i, s[i]=s[i+1]=\cdots=s[j]\rbrace.
\end{equation*}
\label{def_MC}
\end{defn}
Note that the minimal consecutiveness sequence may not be unique. The intuition for concentrating on sequences of minimal consecutiveness is that, under this class of sequences, each target is visited in the shortest time instants compatible with the optimal probability distribution. We next provide a heuristic algorithm with objective to construct a periodic minimal consecutive sequence.
\begin{algorithm}
\caption{Construction of a scheduling sequence of minimal consecutiveness}
\begin{algorithmic}[1]
\REQUIRE{$L$: the sequence length; $n_i$: the number of occurrences of $a_i$ ($i\in \lbrace 1,2,\cdots, N\rbrace$) in a $L$-length sequence.
Without loss of generality, we assume that $n_1\geq n_2 \geq \cdots \geq n_N$.}
\ENSURE{$\lbrace s[k]\rbrace_{k=1}^{L}$}
\STATE{$\lbrace s[k]\rbrace_{k=1}^{n_1}$ $\leftarrow$ Generate a $n_1$-length sequence with all entries as $a_1$.}
\FORALL{$i\in \lbrace 2,\cdots, N\rbrace$}
\STATE{Interpolate an element $a_i$ for every $m_i$'s $a_1$ with equal interpolation interval where $m_i = \lceil{\frac{n_1}{n_i+1}}\rceil$. This can be done by using a backoff counter for each $a_i$. From the beginning of the sequence, the counter is reduced by one whenever an $a_1$ is found and reset to be $\lceil{\frac{n_1}{n_i+1}}\rceil$ as long as $a_i$ is placed.}
\ENDFOR
\end{algorithmic}
\label{alg_seq_MinConsecutive}
\end{algorithm}
\begin{rem}
First of all, the complexity of running Algorithm \ref{alg_seq_MinConsecutive} is $\mathcal{O}(L)$. Next, in our stochastic scheduling strategy, the integer value $n_i$ is obtained from the optimal probability distribution, i.e., $n_i = \lfloor{q_i L}\rfloor$. In order to generate a scheduling sequence precisely matching the distribution $q_i's$ within certain precision, $L$ should be chosen such that $q_iL$ ($i = 1,2,\cdots, N$) is an integer. The rare cases where L is supposed to be infinity, are approximated by taking L large, and a viable sequence is then obtained by periodic continuation. Therefore, besides practical requirements on the length of measuring period, we should choose $L$ by taking into account both the computing capability of the centralized scheduler and the probability matching precision.
\end{rem}
\begin{rem}
Several other relevant deterministic scheduling strategies can be found in the literature, e.g. Round Robin, Pinwheel Scheduling \cite{Chan_pinwheel_92,Chan_pinwheel_93}. However, a careful review finds that they are not suitable/extendable to our specific settings.
\end{rem}
\subsection{A Distributed Scheduling Implementation}
\indent The above algorithm lends itself to a distributed implementation. We assume that the estimators are strongly connected through a network where the communications between
any two estimators are error-free. Under the assumption that each estimator is capable of sensing the channel state, i.e., idle or occupied, we propose a heuristic distributed scheduling mechanism to approximate the result of Algorithm \ref{alg_seq_MinConsecutive}. This distributed scheduling scheme is derived based on the well-known mechanism - carrier sense multiple access with collision avoidance (CSMA/CA). As used in CSMA/CA, our proposed multiple access mechanism relies on backoff timers as well. Specifically, each estimator generates a backoff timer as $T_i = \frac{\alpha}{q_i}$ where $\alpha$ is chosen such that $T_i << \tau$ ($\tau$ refers to the sampling period of DTLTI systems). Then each estimator regularly senses the transmission channel during its backoff timer. If the channel is sensed to be ``idle'' and the timer of estimator $i$ goes off, estimator $i$ begins to use the channel for observing the $i$-th DTLTI system. If the channel is sensed to be busy, then the backoff timer must be frozen until the channel becomes free again. Remark that no synchronization is necessary among backoff timers. With small probability, collision may happen among estimators. That is, backoff timers of two or more estimators go off at the same moment. To deal with this problem, the backoff timer of $i$-th estimator (subject to collision) can be adjusted as $T_i = \frac{\alpha}{q_i-\epsilon_i}$ where $\epsilon_i>0$ is randomly chosen and $\epsilon_i << q_i$. After estimator $i$ uses the channel once, $T_i = \frac{\alpha}{q_i-\epsilon_i}$ is set to be $T_i = \frac{\alpha}{q_i}$ immediately. We remark that, for a large period $T$, the number of observations on $i$-th DTLTI system is approximately $q_i\frac{T}{\alpha}$ where $\frac{T}{\alpha}$ can be treated as a normalizing factor. Thus the optimal probability distribution $q_i's$ is preserved.
\section{Examples and Simulations}
\begin{figure}
\includegraphics[scale=0.63]{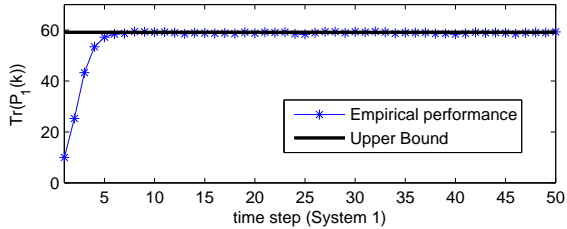}
\caption{Empirical estimate error covariances for DTLTI system (\ref{sys1.eq}). The performance for system $2$ has the same feature. The empirical performance curve is obtained by averaging $5000$ Monte Carlo simulations.}
\label{fig01}
\end{figure}

In this section, we present some simulation results to verify our stochastic scheduling strategy and algorithm.\\
\\
{\em 6.1\; Example A}\\
\\
Consider a single sensor for measuring two DTLTI systems with the following state space representations,
\begin{equation}\label{sys1.eq}
A_{1} = \begin{bmatrix} 0 & 1\\ -0.49 & 1.4  \end{bmatrix}\hspace{-1mm},\, C_{1} = \begin{bmatrix}1\\0\end{bmatrix}^{'}\hspace{-2mm},\,\,Q_{1} = \begin{bmatrix} 5 & 0\\ 0 & 5  \end{bmatrix}\hspace{-1mm},\, R_1=0.5
\end{equation}
\begin{equation}\label{sys2.eq}
A_{2} = \begin{bmatrix} 0 & 1\\ -0.72 & 1.7  \end{bmatrix}\hspace{-1mm},\,C_{2} = \begin{bmatrix}1\\0\end{bmatrix}^{'}\hspace{-2mm},\,\,  Q_{2} = \begin{bmatrix} 1 & 0\\ 0 & 1  \end{bmatrix}\hspace{-1mm}, R_2=1\\
\end{equation}
\begin{table}
\caption{Results of the example}
\begin{tabular}{|c|c|c|c|c|}
  \hline
   & Opt. $q_i$'s & Opt. Cost & Emp. Cost& MC Cost\\
   \hline
  OP & $0.674, 0.326$ & $59.1$ & $58.7$ & $55.7$\\
  \hline
\end{tabular}
\label{res.tab}
\end{table}
The results of applying  Algorithm \ref{alg_bisection_distributed} are reported in Table \ref{res.tab}.
In the table, we also report the empirical cost obtained by pre-assigning the optimal probability distribution to the sensor
and computing the empirical error state covariance of each system by accordingly generating random switching sequences.
This simulation result  verifies that OP provides an effective upper
bound on the expected estimate error covariances of our stochastic scheduling strategies. This is further exemplified in
Fig.\ref{fig01}, which shows the empirical expected estimate error covariance of system (\ref{sys1.eq}).\\

\begin{figure}
\includegraphics[scale=0.45]{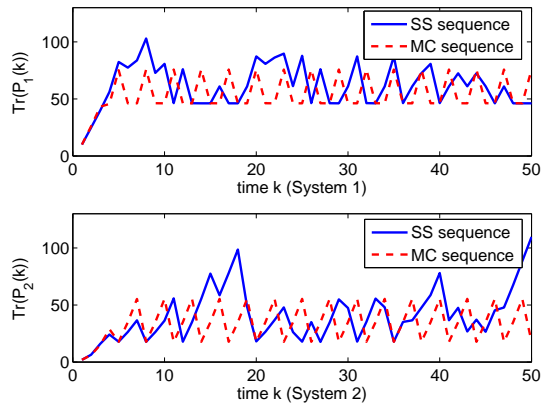}
\caption{Performance of minimal consecutiveness sequence \textit{vs.} stochastic scheduling sequence.}
\label{fig02}
\end{figure}

\begin{figure}
\begin{center}
\includegraphics[scale=.6]{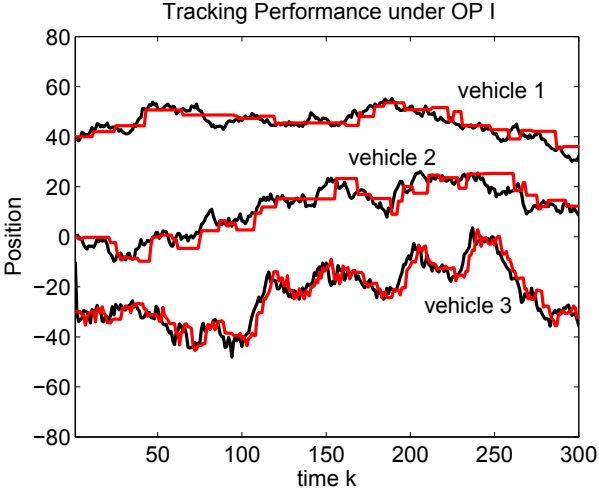}
\end{center}
\caption{Positions of random-walk vehicles and tracking paths.}
\label{figure:tracking}
\end{figure}

{\em 6.1.1\; Minimal consecutiveness (MC) \textit{v.s} stochastic scheduling (SS)} \\
\\
As shown in Fig. \ref{fig02}, the blue solid curves show the evolution of the error covariance
w.r.t time $k$. We see several high peaks as a result of consecutive loss of observations.
The dash lines show that the MC sequence provides much smoother error covariance which is quite necessary in certain applications.
Furthermore, the objective value provided by MC sequence, $55.7$, is smaller than the empirical expected
estimate error covariances, $58.7$, provided by the SS sequences.
This performance enhancement is the result of ``smoothness'' provided by the MC sequence.
In a word, the MC sequence provides a better upper bound than that provided by the SS sequences (i.e. the optimal objective value $\mu_s$ in (\ref{Intro_OP02})). \\
Finally, for completeness, we compare the MC cost with the cost of the (tree search based) sliding window algorithm (see \cite{Gupta04}). It is worth noting that the optimal cost of the original deterministic scheduling problem can be approximated by sliding window algorithm with sufficiently large window size, at the expense of huge computing cost. In this example, the MC cost is $\simeq 55.7$ (in Table \ref{res.tab}), while the sliding window algorithm has the best cost of $\simeq 57.9$ with window size up to $15$.
In this case, we see
that MC not only has better performance than the
sliding window performance, but also can be determined off-line, has low computing complexity and is simple to implement, and scales much better with the number of systems to schedule.\\
\\
{\em 6.2\; Example B}\\
\indent In this example, we consider three random-walk vehicles in an area and a single sensor equipped with a camera is used for tracking their $1-D$ positions. The dynamics of their positions are evolving as (\ref{Model_sec4_01}) with $a_i=1$.
But they are subject to different process noises, measurement noises and delays. Here we assume that these vehicles have $1$, $2$ and $2$ time-step measurement delays, respectively.
Then we have expanded state space systems
\begin{equation*}
\begin{split}
&A_{1} = \begin{bmatrix} 0 & 1\\ 0 & 1  \end{bmatrix},\quad B_{1}=\begin{bmatrix}0\\1 \end{bmatrix},\quad C_{1} = \begin{bmatrix}1\\0\end{bmatrix}^{'}\\
&A_{i} =\begin{bmatrix} 0 & 1 & 0 \\ 0 & 0 & 1 \\ 0 & 0 & 1 \end{bmatrix},\quad B_{i}= \begin{bmatrix}0\\0\\1\end{bmatrix} \quad C_{i}=\begin{bmatrix}1\\0\\0\end{bmatrix}^{'}\; i=2,3\\
&Q_{1}=1 \quad Q_{2}=2 \quad Q_{3}=5, \quad R_{1}=R_{2}=R_{3}=1.\\
\end{split}
\end{equation*}
By running the proposed distributed algorithms, the solutions of OP is $[q_{1},q_{2},q_{3}]=[0.0649, 0.1612, 0.7739]$. We assume the existence of a centralized scheduler which constructs a random scheduling sequence accordingly.
The tracking performances are shown in Fig. \ref{figure:tracking}. The tracking paths (red curves)
are shown in comparison with actual time-varying positions (black curves) of random-walk vehicles.
The flat segments of red curves imply that no measurement is taken in this time slot and the estimator simply
propagates the state estimate of the previous time-step. In Fig. \ref{figure:tracking}, it is shown that the actual path of vehicle $1$ changes slowly, correspondingly the attention given to this sensor is small and the estimate path is updated rarely. In comparison, the tracking paths of vehicle $2$ and vehicle $3$ match the actual positions much better even though the flat segments in tracking path of vehicle $2$ is occasionally visible.
\section{Conclusion}
In this paper, we have proposed a stochastic sensor scheduling problem where the states of $N$ targets (i.e., DTLTI systems) need to be estimated. Due to some practical sensing constraints, at each time instant, only one target can be measured by the sensor. The basic idea of stochastic scheduling is that at each time instant a target is randomly chosen for observation according to a pre-assigned probability distribution. We have relaxed the problem to a convex optimization problem and proposed a distributed algorithm to obtain the optimal probability distribution with the objective to minimize the maximal estimate error among all targets. Then we proposed centralized and distributed scheduling implementation schemes. Finally, we presented simulation results to verify our stochastic strategy.

\section{Appendix}
{\em A. Proof of Theorem \ref{upperBound}}\\
\Proof
\indent Firstly, we define an upper bound $\bar{\mu}_d$ on the performance $\mu_d$ of the deterministic scheduling problem (\ref{Intro_OP01}). For a given $q_1,q_2,\cdots, q_N$ ( $\sum_{i=1}^N q_i=1$), we generate switching sequences with the distribution (\ref{xi}). Clearly, these sequences are not necessarily optimal. Let
\begin{equation}
\begin{split}
\quad &\bar{\mu}_d(q_1,q_2,\cdots, q_N)= \displaystyle\min_{q_i, i=[1..N]}\max_i(\limsup_{T\rightarrow \infty} \frac{1}{T}\sum_{k=1}^{T} Tr(P_i[k]))\\
\end{split}
\label{proof_thm1_equ01}
\end{equation}
with constraints (\ref{errorCov_evolvingEqu}) and (\ref{xi}). Note that given the optimal distribution of $q_1,q_2,\cdots, q_N$ this upper bound still holds since the scheduling sequence is chosen randomly according to (\ref{xi}).
Furthermore, under the assumption, the limit value in (\ref{proof_thm1_equ01}) exists almost surely.
This is because for a stochastic process $\lbrace P_i[k]\rbrace$ ($P_i[0]=\Sigma$) satisfying (\ref{errorCov_evolvingEqu}), there exists an ergodic stationary process $\lbrace \bar{P}_i[k]\rbrace$ satisfying (\ref{errorCov_evolvingEqu})
with
\begin{equation}
\lim_{k\rightarrow \infty} ||P_i[k]-\bar{P}_i[k]|| = 0. \quad \textit{a.s.}
\end{equation}
See \cite{Bougerol} Theorem 3.4.
That is,
\begin{equation}
\begin{split}
\lim_{T\rightarrow \infty} \frac{1}{T}\sum_{k=1}^{T} Tr(P_i[k]) &= \lim_{T\rightarrow \infty} \frac{1}{T}\sum_{k=1}^{T} Tr(\bar{P}_i[k]). \quad \textit{a.s.}\\
\end{split}
\end{equation}

Thus, (\ref{proof_thm1_equ01}) can be rewritten as
\begin{equation}
\begin{split}
\quad &\bar{\mu}_d(q_1,q_2,\cdots, q_N)= \displaystyle\min_{q_i, i=[1..N]}\max_i(\lim_{T\rightarrow \infty} \frac{1}{T}\sum_{k=1}^{T} Tr(P_i[k]))\\
\end{split}.
\end{equation}
Next, for every sensor, its selection process can be viewed as a binary tree where the variance of observation noise on one branch is infinity if the sensor is not selected at this time instant. This model is a special case of the one (tree topology) considered in \cite{Mo_stochastic11}. Under the current assumptions it is straightforward to check that all assumptions in Theorem $1$ of \cite{Mo_stochastic11} hold. As a result, inequality (10) in \cite{Mo_stochastic11} holds in our case, i.e.,
\begin{equation*}
\lim_{T\rightarrow \infty} \frac{1}{T}\sum_{k=1}^{T}Tr(P_i[k])\leq Tr(\lim_{k\rightarrow \infty} \mathbb{E}_{\xi_{i,k}}[P_i[k]]). \quad \textit{a.s.}
\end{equation*}
Recalling that $
\mu_s =\displaystyle\min_{q_i, i=[1..N]} \max_i Tr( \lim_{k\rightarrow \infty} \mathbb{E}_{\xi_{i,k}}[P_i[k]])$,
it follows that $\mu_d \leq \bar{\mu}_d \leq \mu_s$. Namely, the performance $\mu_d$ of the deterministic scheduling problem (\ref{Intro_OP01}) is upper bounded by the performance $\mu_s$ of the stochastic scheduling problem (\ref{Intro_OP02}) \qed\\

\bibliographystyle{IEEEtran}
\bibliography{ref}

\end{document}